\documentclass[12pt]{article}

\def\qed{\quad \vrule height7.5pt width4.17pt depth0pt}

\newtheorem{assumption}{Assumption}
\newtheorem{definition}{Definition}
\newtheorem{theorem}{Theorem}
\newtheorem{lemma}{Lemma}

\begin{document}

\title{The heat equation with multiplicative stable L\'evy noise}
\author{
Carl Mueller$^{1}$
\\Dept. of Mathematics
\\University of Rochester
\\Rochester, NY  14627
\\E-mail:  cmlr@math.rochester.edu
\and
Leonid Mytnik$^{2}$
\\Faculty of Industrial Engineering and Management
\\Technion -- Israel Institute of Technology
\\Haifa 32000, Israel
\\E-mail:  leonid@ie.technion.ac.il
\and
Aurel Stan$^{1}$
\\Dept. of Mathematics
\\University of Rochester
\\Rochester, NY  14627
\\E-mail:  astan@math.rochester.edu}
\date{}
\maketitle

\footnotetext[1]{Supported by an NSF grant.

$\;\mbox{}^{2}$\hspace*{-1.3mm} Research supported in part by the Israel Science Foundation
(grant No.~116/01-10.0).


AMS 1991 {\em subject classifications}
Primary, 60H15; Secondary, 35R60, 35K05.}

\newpage

\abstract{ We study the heat equation with a random potential term.  The 
potential is a one-sided stable noise, with positive jumps, which does not 
depend on time.  To avoid 
singularities, we define the equation in terms of a construction similar to 
the Skorokhod integral or Wick product.  We give a criterion for existence 
based on the dimension of the space variable, and the parameter $p$ of the 
stable noise.  Our arguments are different for $p<1$ and $p\geq1$.}

\newpage

\section{Introduction}
\setcounter{equation}{0}

Our goal is to study the heat equation with stable L\'evy noise $\dot{L}(x)$ 
which depends only on space.
\begin{assumption}
Throughout the paper, we assume that $\dot{L}$ has only positive jumps.  
\end{assumption}
Our noise is multiplicative, in the sense that 
it is multiplied by $u$ in the equation.  
\begin{eqnarray}
\label{1.1}
\frac{\partial u}{\partial t} &=& \Delta u + u\diamond\dot{L}(x)  \\
u(0,x) &=& u_0(x).  \nonumber
\end{eqnarray}
Here, $\dot{L}(x)$ is a one-sided stable noise of index $p\geq1$ defined on 
$x\in\mathbf{R}^d$.  The product 
$u\diamond\dot{L}(x)$ is related to the Skorokhod integral or Wick product.  
We will give precise conditions later.  
Our main goal is to 
give conditions on $p$ and the dimension $d$ such that the solution exists and 
is unique.  In fact, for $p<1$, the stochastic integral is defined in a 
different way and can include more singular terms.  
In this case, we use the notation
\begin{eqnarray}
\label{1.2}
\frac{\partial u}{\partial t} &=& \Delta u + u\star\dot{L}(x)  \\
u(0,x) &=& u_0(x).  \nonumber
\end{eqnarray}

There are many papers concerning the heat equation with a multiplicative noise 
term (see \cite{par93}).  In most cases, however, the noise depends on time, 
and is even independent from one time point to another.  This is a very 
natural assumption from the point of view of It\^o integration, since the 
relevant integrals often yield martingales.  For example, the equation 
\begin{eqnarray}
\label{1.3}
\frac{\partial u}{\partial t} &=& \Delta u + f(u)\dot{W}(t,x)  \\
u(0,x) &=& u_0(x)  \nonumber
\end{eqnarray}
where $\dot{W}(t,x)$ is space-time white noise,
is usually formulated in terms of an integral equation 
\begin{eqnarray*}
u(t,x) &=& \int_{\mathbf{R}^d}G(t,x-y)u_0(y)dy \\
&&+\int_{0}^{t}\int_{\mathbf{R}^d}G(t-s,x-y)f(u(s,y))W(dyds)
\end{eqnarray*}
where $G(t,x)$ is the heat kernel, and where the final integral can be defined 
in much the same way as the It\^o integral.  

But if the noise does not depend on time, such an integral would be 
anticipating.  Further evidence of the difficulty of noise depending only on 
space is given by the following equation:  
\begin{eqnarray}
\label{1.4}
\frac{\partial u}{\partial t} &=& \Delta u + \delta(x)u  \\ 
u(0,x) &=& u_0(x)  \nonumber
\end{eqnarray}
where $\delta(x)$ is the Dirac delta function.  When $u$ does not depend on $t$, 
this equation has been extensively studied by mathematical physicists.  See 
the book \cite{aghh88} by Albeverio et al., for example.  This equation 
(without $t$ dependence) models the quantum mechanics of a particle under the 
influence of a point interaction.  As usual, we can think of this equation in 
terms of the density of Brownian particles.  The $\delta(x)$ term means that 
Brownian particles at 0 would give birth to many new particles.  But these new 
particles would themselves be close to 0, and so they might give birth to 
other particles, ad infinitum.  Because of this unstable behavior, when the 
dimension $d$ is 2 or higher we need to take an approximate $\delta(x)$ 
function, multiply it by a very small term, and take the limit of the 
approximation.  The case of $x\in\mathbf{R}$ is better behaved.  

Yaozhong Hu has written two papers on time-independent noise, 
\cite{hu01} and \cite{hu02}.  He deals with the equation
\begin{eqnarray*}
\frac{\partial u}{\partial t} &=& \Delta u + u\diamond\dot{W}(x)  \\
u(0,x) &=& u_0(x)
\end{eqnarray*}
where $\diamond$ denotes Skorokhod integration.  One can also think of this 
product as the Wick product, see \cite{jan97} Chapter 3.  $\dot{W}(x)$ is 
white noise in \cite{hu02}, but in \cite{hu01} it is colored Gaussian noise.  
By the Skorokhod integral, roughly speaking, we mean that in the relevant 
stochastic integrals, we should drop diagonal terms such as 
$(\dot{W}(x))^2$.  These terms correspond to repeated influence of $\dot{W}(x)$ 
for the same point $x$, so they are related to the singular behavior in the 
equation (\ref{1.4}).

The goal of this article is to study (\ref{1.1}) for time-independent 
L\'evy noise $\dot{L}(x)$.  The product $\diamond$, like the Wick product, 
involves the deletion of diagonal terms in the appropriate integrals.  But it 
is not defined in the same way as the Wick product for Brownian functionals, 
at least with the usual definitions.  Once again, the reason for using such a 
product is to avoid singular terms such as those coming from (\ref{1.4}).
The product $\star$ used in (\ref{1.2}) is similar, but we delete 
fewer terms.  

Our idea is that, once we understand the equation without the singular terms, 
we can try gingerly adding back the singular terms.  But the understanding of 
the simpler situation must come first.  

In the time dependent situation, there have been several papers involving 
L\'evy noise, such as Kallianpur, Xiong, Hardy, and Ramasubramanian 
\cite{kxhr94}.  In most cases, the noise has been 
additive.  That is, it appears without any multiple of $u$ or $f(u)$.  

In the case $p<1$, and for noise depending on both time and space, 
\cite{mue98} dealt with the equation
\begin{eqnarray}
\label{1.5}
\frac{\partial u}{\partial t} &=& \Delta u + u^\gamma\dot{L}(t,x)  \\
u(0,x) &=& u_0(x).  \nonumber
\end{eqnarray}
If $x\in D$ for a smooth and bounded domain $D\subset\mathbf{R}^d$, and with 
some conditions on $u_0$, short-time existence was obtained in the case
\[
d<\frac{(1-p)\alpha}{\gamma p-(1-p)}.
\]

For $p\geq1$, Mytnik \cite{myt02} found that if
\begin{eqnarray*}
p\gamma &<& \frac{2}{d} + 1,  \\
1 \ <\  p &<&  \min\left(2,\frac{2}{d}+1\right)
\end{eqnarray*}
then there exists a weak solution to (\ref{1.5}).

Both of these papers were motivated by the superprocess with stable branching,
see \cite{daw93}.  We will not describe this superprocess, but it is believed 
that in some weak sense, it is associated with the equation (\ref{1.5}) 
with $\gamma=1/p$.  This motivation led us to study equations with one-sided 
noise.  We leave possible generalizations to the reader.

Our results are also divided into the two cases $p<1$ and $p\geq1$.  Our main 
tools are multiple stochastic integrals, and this is the reason that we 
restrict ourselves to noise terms of the form $f(u)\diamond\dot{L}(x)$ or 
$f(u)\star\dot{L}(x)$ with $f(u)=u$.  
Hu restricted himself to similar equations in the 
Gaussian case, for the same reason.  Of course, multiple stochastic integrals 
have been used many times to study equations such as (\ref{1.3}).  For 
example, see Nualart and Zakai \cite{nz89} and Nualart and Rozovskii 
\cite{nr97}.  But in these articles, the noise always depends on time.  

\section{Theorems}
\setcounter{equation}{0}

In this section we list our main results.  They depend on some definitions 
which will be explained later.  We assume that $u_0(x)$ is a 
bounded function on $D$.  Since our equation is linear, we may assume 
without loss of generality that $\|u_0\|_\infty\leq1$.  

Here are our main theorems.  

\begin{theorem} \label{t1} Suppose that $0<p<1$, and let $u_0(\cdot)$ be a 
bounded function on $D$.  Assume that either $d\leq4$ or
\[
p < \frac{1}{2} + \frac{1}{d-2}, \qquad d\geq 5.
\]
Then (\ref{1.2}) has a unique solution $u(t,x)$.  
\end{theorem}
We are not sure that the conditions on $p,d$ in the above theorem are optimal.  
\begin{theorem} 
\label{t2} 
Suppose that $p\geq1$ and $u_0(\cdot)$ is a 
bounded function on $D$.  
Assume that
\[
p < 1 + \frac{2}{d}.
\]
Then equation (\ref{1.1}) has a unique solution $u(t, x)$.  Moreover, for 
all $t > 0$ and $x \in D$, $u(t, x) \in \mathbf{L}_{loc}^{q}(\Omega)$, for all 
$q \in \left[1, 1 + (2/d)\right)$. \end{theorem}

Here $\mathbf{L}_{loc}^{q}(\Omega)$ denotes the 
space of all random variables $f: \Omega \to \mathbf{R}$, such that, for all 
$\epsilon > 0$, there exists an event $A \subset \Omega$, of probability less 
than $\epsilon$, such that $E[|f1_{A^c}|^{q}] < \infty$.

\textbf{Remarks:}  
\begin{enumerate}
\item In Theorem \ref{t2}, in the Brownian case of $p=2$, the critical 
dimension is $d=2$, which agrees with \cite{hu02}.
\item For $p<1$, if we use the product $\diamond$ instead of $\star$, the 
proof of Theorem \ref{t1} is trivial.  The reader can check that by removing
the large atoms of the measure $\dot{L}(x)$, we can take $q=1$ and use the 
proof of Theorem \ref{t2}.
\end{enumerate}

\section{Basic Definitions}
\setcounter{equation}{0}

We will index sequences of real numbers, such as $y_i$, by subscripts.  A
sequence of vectors in $\mathbf{R}^n$ will be indexed by superscripts.  Thus,
$x^{(i)}$ is the i$^{\rm th}$ vector in the sequence, and $x^{(i)}_k$ is the
k$^{\rm th}$ component of the i$^{\rm th}$ vector.

If it is not stated otherwise, $C$ in this paper will denote a constant which value my change from line to line. 

Throughout the paper, $G(t,x)$ denotes the fundamental solution of the heat 
equation for $x\in\mathbf{R}^d$,
\begin{eqnarray*}
\frac{\partial u}{\partial t} &=& \Delta u  \\
u(0,x) &=& \delta(x).
\end{eqnarray*}
Explicitly,
\[
G(t,x) = \frac{1}{(4\pi t)^{d/2}}\exp\left(-\frac{x^2}{4t}\right).
\]
Given a domain $D\subset\mathbf{R}^d$, we let $G(t,x,y)$ be the fundamental 
solution of the heat equation in $D$, with Dirichlet boundary conditions.  
That is, $G(t,x,y)$ satisfies 
\begin{eqnarray*}
\frac{\partial u}{\partial t} &=& \Delta u \\
u(0,x) &=& \delta_y(x)  \\
u(t,x) &=& 0 \hspace{.5in} \mbox{for $x\in\partial D$}.
\end{eqnarray*}
We note the elementary inequality 
\[
0 \leq G(t,x,y) \leq G(t,x-y).
\]
In other words, when $G$ has two arguments, it is the heat kernel on 
$\mathbf{R}^d$, and when it has three arguments, it is the heat kernel on $D$.  

Finally, we write 
\[
\mathbf{T}_n = \mathbf{T}_n(t) 
:= \{(s_1,\dots,s_n): 0<s_1<\dots<s_n<t\}.
\]

\section{The Noise}
\setcounter{equation}{0}
\label{the-noise}

Roughly speaking, our noise $\dot{L}(x)$ will be the analogue of $dL(t)/dt$ 
where $L(t)$ is a one-sided stable process of index $p\in(0,2)$, and 
\newline $t\in[0,\infty)$.  Of course, 
this derivative only exists in the generalized sense.  For example, we could 
use the theory of Schwartz distributions.  

We will start with a review of the one-sided stable L\'evy processes, and 
refer the reader to Bertoin's book \cite{ber96} for details.  There is a 
basic difference between the cases $p<1$ and $p\geq 1$.  

For $p<1$, we can construct $L(t)$ via a Poisson process, as follows.  
Consider the following measure on $(t,y)\in[0,\infty)^2$ given by
\[
\nu(dydt) = c_p y^{-(p+1)}dydt,
\]
where $c_p$ is a constant depending on $p$.  This is the same constant that 
occurs in the L\'evy measure of the one-sided stable processes 
(see \cite{ber96}).
Let $Y(\cdot)$ be a Poisson random measure on $[0,\infty)^2$ 
with intensity $\nu$.  Then, the $p$-stable process $L(t)$ is defined as 
\[
L_t = \int_{[0,t]}\int_{0}^{\infty}yY(dyds).
\]
Note that if we let $(t_i,y_i)$ be the locations of the points (or atoms) of 
the Poisson measure $Y$, then
\[
L_t = \sum_{t_i\leq t}y_i.  
\]
In a similar way, we define the random measure 
$L(A): A\subset[0,\infty)$ 
as 
\[
L(A) = \sum_{t_i\in A}y_i
\]
and observe that
\[
L_t = \int_{0}^{t}L(ds).
\]

For $p\geq1$, we must introduce an approximation and compensation procedure, 
as follows.  Let 
\[
Y_n(dydt) = Y(dydt)\mathbf{1}(y\geq 1/n).
\]
Note that if we define
\[
\nu_n(dydt) := \nu(dydt)\mathbf{1}(y\geq 1/n)
\]
then $Y_n$ becomes a Poisson measure with intensity $\nu_n$.  Note that 
\[
\nu_n([0,t]\times[0,\infty)) = 
\nu([0,t]\times[1/n,\infty)) = \frac{c_p}{p} n^p t.  
\]
Next, let
\begin{eqnarray*}
L^{(n)}_t &=& \int_{0}^{t}\int_{0}^{\infty}yY_n(dyds) - 
t\int_{0}^{\infty}yd\nu_n(y)  \\
&=& \int_{0}^{t}\int_{0}^{\infty}yY_n(dyds) - \frac{c_p}{p-1}n^{p-1}t.  
\end{eqnarray*}
There is a problem with $p=1$.  We leave the details to the reader, 
but roughly speaking, we would truncate the measure by removing large jumps.  

Note that $L_t^{(n)}$ is a martingale with respect to its own filtration.  The 
$p$-stable process $L(t)$ is the limit of the processes $L_t^{(n)}$ as 
$n\to\infty$.  For details of this limit, we refer the reader to \cite{ber96}.
We could define the random measure $L(dt)$ in a similar way.  

The notation $\dot{L}(t)$ stands for the the density of the measure $L$, even 
though this is only defined in the generalized sense.  This is similar to the 
convention of writing $\dot{B}(t)$ for the derivative of Brownian motion.  

Our construction of $L(A): A\subset\mathbf{R}^d$ is similar.  First consider 
the case $p<1$.  Fix a bounded open region $D\subset\mathbf{R}^d$;  we will 
give some conditions on $D$ later.  For $(x,y)\in B\subset D\times[0,\infty)$, 
let $Y(B)$ be the Poisson measure with intensity $dx\nu(dy)$.  Let 
$(x^{(i)},y_i)$ denote the locations of the points (atoms) of the Poisson 
measure $Y$.  For $A\subset D$, we define
\[
L(A) = \sum_{x^{(i)}\in A}y_i.  
\]

For $p\geq 1$, we must approximate.  Given $n$, let $Y,\nu$ be defined as above,
and let
\begin{eqnarray*}
Y_n(dxdy) &=& \mathbf{1}(y\geq 1/n)Y(dxdy)   \\
dx\nu(dy) &=& c_p y^{-(p+1)}dxdy       \\
dx\nu_n(dy) &=& \mathbf{1}(y\geq 1/n)dx\nu(dy).
\end{eqnarray*}
As above, for $A\subset D$, we let
\begin{eqnarray*}
L_n(A) &=& 
\sum_{x^{(i)}\in A}y_i \mathbf{1}(y_i\geq 1/n) 
-\int_{0}^{\infty}\int_{A}ydx\nu_n(dy)  \\
&=& \sum_{x^{(i)}\in A}y_i \mathbf{1}(y_i\geq 1/n) 
- m(A)\frac{c_p}{p-1}n^{p-1}
\end{eqnarray*}
where $m(A)$ denotes the Lebesgue measure of $A$.  Again, we leave the case 
$p=1$ to the reader.  To get $L(dx)$, we 
again take a limit as $n\to\infty$.  Since this procedure is not very different 
than for the one-sided stable processes with $p\in[1,2)$, we refer the reader 
to \cite{ber96}.

For the future, we label the atoms of the measure $L$ by 
$(x^{(i)},y_i): x^{(i)}\in\mathbf{R}^d$.  Then, $x^{(i)}$ is the location of 
the atom, and $y_i$ is the mass of the measure $L(\{x^{(i)}\})$.

Fix $K>0$.  Since $D$ is a bounded domain, there are only a finite number of 
atoms $(x^{(i)},y_i)$ with $x^{(i)}\in D$ and $y_i>K$.  Let $L_K$ denote $L$ 
with the preceding atoms deleted.  That is,
\[
L_K = L - \sum_{y_i>K}y_i\delta_{x^{(i)}}.
\]
\begin{assumption}
In the succeeding sections, we will tacitly assume that $\dot{L}$ is 
actually $\dot{L}_K$.  Furthermore, let $\mathcal{A}_K$ denote the event that 
$L=L_K$, that is, there are no atoms larger than $K$.  
\end{assumption}
Note that
\[
\lim_{K\to\infty}P\left(\mathcal{A}_K\right) = 
\lim_{K\to\infty}P\left(L = L_K\right) = 1
\]
 From now on, we will replace $L$ by $L_K$, and let $K\to\infty$ at the end.

\section{Multiple Stochastic Integrals}
\setcounter{equation}{0}
\label{precise}

In the one-dimensional case, such integrals have been considered earlier, for 
example in \cite{rw86}.
We suppose that $D\subset\mathbf{R}^d$ is an open, bounded domain with smooth 
boundary.  As described earlier, for $p\in(0,2)$, let $L(A)$ be the 
one-sided stable random measure of index $p$, defined for Borel sets 
$A\subset D$.  

Assume that $p\geq1$.  Now we describe our multiple integrals with respect to 
the noise $\dot{L}$, recalling that $\dot{L}(x)dx$ is just another notation 
for $L(dx)$.  Assume that we have a  symmetric function 
$f_n(x^{(1)},\dots,x^{(n)})$ defined for $x^{(1)},\dots,x^{(n)}\in D$.  
Define $\widetilde{D}_{n} \subset D^n$ as follows.  
\[
\widetilde{D}_n = \{(x^{(1)}, x^{(2)}, \dots, x^{(n)}) \in D^{n} ~:~  
x_{1}^{(1)} < x_{1}^{(2)} < \dots < x_{1}^{(n)}\}.
\]
We regard $x_1$ as the time variable, so that the following integral is 
defined in the It\^o sense.  
\[
I_n(f_n) = n!\int_{\widetilde{D}_{n}}f_n(x^{(1)},\dots,x^{(n)})
L(dx^{(1)})\dots L(dx^{(n)}).
\]
Sufficient conditions for the existence of the integral $I_n(f_n)$ are given 
in Lemma \ref{poissonian_Wiener_Ito}.

Since $f_n$ is symmetric, $I_n(f_n)$ is invariant under permutations of the 
indices $i=1,\dots,n$ of the $x^{(i)}$.  
For $f=(f_0,f_1,\dots)$, where $f_n$ is a symmetric function of the variables
$x^{(1)},\dots,x^{(n)}$, let 
\[
I(f) = \sum_{n=0}^{\infty}I_n(f_n).
\]
Here, $f_0$ is interpreted as a constant, and $I_0(f_0)=f_0$.  

For future use, we define the symmetrization of a function as follows.  
\[
{\rm sym}(f)(x^{(1)},\dots,x^{(n)}) 
= \frac{1}{n!}\sum_{\sigma}f(x^{(\sigma(1))},\dots,x^{(\sigma(n))})
\]
where the sum is taken over all permutations $\sigma$ of the indices 
$1,\dots,n$.  If we wish to symmetrize over only some of the variables, we 
write those in a subscript.  
\begin{eqnarray*}
\lefteqn{ {\rm sym_{x^{(1)},\dots,x^{(n)}}}(f)(x^{(1)},
\dots,x^{(n)},z^{(1)},\dots,z^{(m)}) }\\
&&= \frac{1}{n!}\sum_{\sigma}f(x^{(\sigma(1))},\dots,x^{(\sigma(n))},
z^{(1)},\dots,z^{(m)})
\end{eqnarray*}

Note that the integral $I(f)$ has the following property.  Because of the 
ordering we have imposed on the variables $x^{(i)}$, it follows that $I(f)$ 
depends linearly on the masses $y_i$ of the atoms of the measure $L$.  
That is, $I(f)$ can have terms like $y_1y_2$, but no quadratic terms like 
$y_1^2$.

We note in passing that our multiple stochastic integrals are related to the 
Wick product for the compensated Poisson process $Y$;  see Lytvynov 
\cite{lyt03}, for example.   

For the case $p<1$, we propose a different integral $J(f)$, which allows for 
more singularities.  Let $\mathcal{X}_n(D)\subset D^n$ be the set of 
coordinates $(x^{(1)},\dots,x^{(n)})$ with no adjacent indices equal.  In 
other words, $x^{(i)}\ne x^{(i+1)}$ for $i=1,\dots,n-1$.  Let
\[
J_n(f_n) = \int_{\mathcal{X}_n(D)}f_n(x^{(1)},\dots,x^{(n)})
L(dx^{(1)})\dots L(dx^{(n)}).
\]
Since $L(dx^{(1)})\dots L(dx^{(n)})$ is a measure, the integral $J_n(f_n)$ can 
be defined for each $\omega$, provided $f_n$ is integrable with respect to this 
product measure.   

We define $J(f)$ in a similar way:
\[
J(f) = \sum_{n=0}^{\infty}J_n(f_n).
\]

\section{A Framework for the Equation}
\setcounter{equation}{0}

The goal of this section is to set up a framework for discussing the equation, 
so that existence and uniqueness can be rigorously discussed.  In the white 
noise case, there are some related concepts in \cite{nr97} and \cite{mr98}.

We define 
\begin{equation}
{\cal I}_n \equiv \{ I_n(f_n) : f_n \ \mbox{is a symmetric function on} \; 
D^n\}.
\end{equation} 

\begin{definition}
Let $g(\cdot,\cdot)$ be a measurable function  on $D\times \Omega$, such that 
   \mbox{$g(x,\cdot)\in {\cal I}_n$} for almost every $x\in D$. That is, there exists a function $f_n$ on $D^{n+1}$ such that 

\[
g(x,\cdot)=n!\int_{\widetilde{D}_n}f_n(x^{(1)},\dots, x^{(n)},x)
L(dx^{(1)})\dots L(dx^{(n)}),\;\; {\rm a.e.}\; x,
\] 
and 
$ f_n (\cdot,x)$ is a  symmetric function on $D^n$ for almost every $x$.

Then we define
\begin{eqnarray*}
\lefteqn{\int_D g(x)\diamond L(dx)
 }
\\
&\equiv & (n+1)! \int_{\widetilde{D}_{n+1}}
{\rm sym}_{x^{(1)},\dots,x^{(n+1)}}
    f_n(x^{(1)},\dots,x^{(n)}, x^{(n+1)})
\\
&&\mbox{}\hspace{2cm} 
L(dx^{(1)})\dots L(dx^{(n+1)}).  
\end{eqnarray*} 
\end{definition}
Set 
\[
{\cal I}=\left\{ \sum_{n=0}^{\infty} g_n : g_n\in {\cal I}_n\;  {\rm and \; the\; sum \; is \; converging \; in \; probability}\right\}.
\] 

\begin{definition}
Let $g(\cdot,\cdot)$ be a measurable function  on $D\times \Omega$, such that 
\mbox{$g(x,\cdot)\in {\cal I}$} for almost every $x\in D$.  That is, there 
exists a sequence $\{g_n\}_{n\geq0}$ such that 
\[
g_n(x,\cdot)\in {\cal I}_n\,, \;\; n\geq 0, \; {\rm a.e.}\; x\in D,
\] 
and 
\[
g(x,\cdot)=
 \sum_{n=0}^{\infty} g_n(x,\cdot),\;\;  {\rm a.e.} \; x\in D.
\] 
Then we define
\[
\int_D g(x)\diamond L(dx)\equiv \sum_{n=0}^{\infty} \int_D g_n(x)\diamond L(dx),
\]
where  the  sum is converging  in  probability.
\end{definition}

Next definition is analogous to the definition of Hu, from \cite{hu01}. 

\begin{definition}
A measurable function  
$u: \mathbf{R}_{+}\times D\times \Omega \mapsto \mathbf{R}$ is 
called a solution to the equation (1.1) if $u(t,x)\in {\cal I}$ for every 
$x\in \mathbf{R}^d$ and $t\in (0,T]$ and the following equation is satisfied:
\begin{eqnarray}
\label{6.1}
u(t,x) &=& \int_{D}G(t,x,y)u_0(y)dy   \\
&&+ \int_{D}\left(\int_{0}^{t}G(t-s,x,y)u(s,y)ds\right)\diamond L(dy),\nonumber\\
&&\qquad \forall  x\in \mathbf{R}^d,
 t\in (0,T].  \nonumber
\end{eqnarray}
\end{definition}

Our strategy is to expand $u(t,x)$ in a recursively defined series:  
\[
u(t,x) = \sum_{n=0}^{\infty}u_n(t,x),
\]
where
\begin{eqnarray}
\label{6.2}
u_0(t,x) &=& \int_{D}G(t,x,y)u_0(y)dy  \\
u_{n+1}(t,x)
&=& \int_{0}^{t}\int_{D}G(t-s,x,y)u_n(s,y)\diamond L(dy) ds.    \nonumber
\end{eqnarray}
Note that $u_0$ has two meanings.  If $u_0(x)$ is a function of one variable, 
then it is the initial value for our SPDE (\ref{1.1}).  If $u_0(t,x)$ 
depends on two variables, then it is the first term of the expansion  for
$u(t,x)$.  

To deal with the case $p<1$, we define the operation $\star$.  
We define 
\begin{equation}
{\cal J}_n \equiv \{ J_n(f_n) : f_n \ \mbox{is a function on} \; D^n\}.
\end{equation} 

\begin{definition}
Let $g(\cdot,\cdot)$ be a measurable function  on $D\times \Omega$, such that 
   \mbox{$g(x,\cdot)\in {\cal J}_n$} for almost every $x\in D$. That is, there exists a function $f_n$ on $D^{n+1}$ such that 

\[
g(x,\cdot)=\int_{\mathcal{X}_n(D)}f_n(x^{(1)},\dots, x^{(n)},x)
L(dx^{(1)})\dots L(dx^{(n)}),\;\; {\rm a.e.}\; x,
\] 
and 
$ f_n$ is a function on $D^{n+1}$.  

Then we define
\begin{eqnarray*}
\lefteqn{\int_D g(x)\star L(dx)
 }
\\
&\equiv & \int_{\mathcal{X}_{n+1}(D)}
    f_n(x^{(1)},\dots, x^{(n+1)})
L(dx^{(1)})\dots L(dx^{(n+1)}).  \nonumber
\end{eqnarray*} 
\end{definition}

The other definitions for $p<1$ are completely analogous to the case $p\geq1$, 
with $\star$ instead of $\diamond$.  We leave this part to the reader.

\section{Proof of the existence part of \\ \mbox{}\hspace{0.8cm} Theorem \ref{t1} ($p<1$)}
\setcounter{equation}{0}

\subsection{Calculations for the heat equation}

As mentioned at the end of the Section \ref{the-noise}, we will assume that 
$\mathcal{A}_K$ occurs, so that $\dot{L}=\dot{L}_K$.  We suppose that $K=2^N$ 
for some integer $N$.  We will obtain conclusions that hold almost surely for 
each integer $K=2^N$.  Since $P(\mathcal{A}_K)\to1$ as $K\to\infty$, 
our assertion will almost surely hold for $L$.  

We will use the notation of Mueller~\cite{mue98}.  We call the atom at 
$x^{(i)}$ with mass $y_i$ a particle of type $n$ if 
\[
2^{-(n+1)} < y_i \le 2^{-n}.
\]
Since there are no atoms of mass greater than $K=2^{N}$, we need only consider 
particles of type $n\geq-N$.  

We will use the expansion defined in (\ref{6.2}), with $\star$ instead 
of $\diamond$.  Then $u_0(t,x)$ is the solution of the heat equation with 
initial function $u_0(x)$.  

\subsection{Distance between points}

We start with an elementary lemma about Poisson probabilities.
\begin{lemma}
\label{poisson-est}
Let $X$ be a Poisson random variable with parameter $\lambda$.  
Then for any nonnegative integer $n_0$ we have
\[
P(X\geq n_0) \leq \frac{\lambda^{n_0}}{n_0!}.
\]
As a consequence
\[
P(X\geq n_0) \leq \frac{\lambda^{n_0}}{n_0^{n_0}}e^{n_0}.
\]
\end{lemma}

\noindent
\textbf{Proof.}
\begin{eqnarray*}
P(X\geq n_0) &=& \sum_{n=n_0}^{\infty}P(X=n)  \\
&=& \sum_{n=n_0}^{\infty}\frac{\lambda^n}{n!}e^{-\lambda}  \\
&=& \frac{\lambda^{n_0}}{n_0!}\sum_{n=n_0}^{\infty}
   \frac{\lambda^{n-n_0}}{(n-n_0)!}\cdot\frac{n_0!(n-n_0)!}{n!}e^{-\lambda}  \\
\end{eqnarray*}
But for all $n\geq n_0$, we have
\[
\frac{n_0!(n-n_0)!}{n!} = \frac{1}{{n\choose n_0}} \leq 1.
\]
Thus,
\begin{eqnarray*}
P(X\geq n_0) &\leq&
   \frac{\lambda^{n_0}}{n_0!}\sum_{n=n_0}^{\infty}
   \frac{\lambda^{n-n_0}}{(n-n_0)!}e^{-\lambda}  \\
&=&
   \frac{\lambda^{n_0}}{n_0!}\sum_{k=0}^{\infty}
   \frac{\lambda^{k}}{k!}e^{-\lambda}  \\
&=& \frac{\lambda^{n_0}}{n_0!}.
\end{eqnarray*}
Since 
\begin{eqnarray*}
e^{n_0} &=& 1 + \frac{n_0}{1!} + \frac{n_0^2}{2!} + \dots + 
         \frac{n_0^{n_0}}{n_0!} + \dots    \\
&\geq& \frac{n_0^{n_0}}{n_0!}
\end{eqnarray*}
we have
\[
P(X\geq n_0) \leq \frac{\lambda^{n_0}}{n_0!} 
\leq \lambda^{n_0}\frac{e^{n_0}}{n_0^{n_0}}.
\]

\qed

Suppose that $D\subset B$, where $B$ is a cube in $\mathbf{R}^d$ with side
$r\geq1$.  Fix $a,\ell>0$.  Consider the set of closed cubes of side length 
$\ell$ whose centers lie in the scaled lattice $a\mathbf{Z}^d$.  Let 
$\mathbf{C}_\ell(B,a)$ be the subset of such cubes which intersect $B$.  

Suppose that $x,y\in D$ is a pair of points whose distance $|x-y|\le a$.  We 
claim that there exists a cube in $\mathbf{C}_{3a}(B,a)$ which contains both 
$x$ and $y$.  Indeed, since the cubes $\mathbf{C}_a(B,a)$ cover $B$, it 
follows that $x\in S_a$ for some cube $S_a\in\mathbf{C}_a(B,a)$.  Let $S_{3a}$ 
be the cube with the same center as $S_a$, but with side length $3a$.  It is 
easy to see that every point within distance $a$ of $S_a$ lies in $S_{3a}$.  
Therefore, $S_{3a}$ is a cube from the set $\mathbf{C}_{3a}(B,a)$ which 
contains both $x$ and $y$.  

Furthermore, let $R\subset\mathbf{R}^d$ and let $N(R,n)$ be 
the number of particles of type $n$ which lie in $R$.  
Let $|R|$ be the Lebesgue measure of $R$.  Note that $N(R,n)$ is a Poisson 
random variable with parameter 
\begin{eqnarray}
\label{7.1}
\lambda(R,n) &=& |R|\int_{2^{-(n+1)}}^{2^{-n}} \nu(dy)     \\
&=& |R|\int_{2^{-(n+1)}}^{2^{-n}} \frac{c_p}{y^{p+1}} dy    \nonumber\\
&=& C_p |R| 2^{np}.  \nonumber
\end{eqnarray}
where $C_p=c_p(2^p-1)/p$.  

For $a,M>0$, let $A(n,m,M,S)$ be the event that within $S\cap B$, there is at 
least one particle of type $n$ and at least $M$ particles of type $m$, where 
$S$ is a cube in $C_{3a}(B,a)$.  
Let
\[
\tilde{S} = S\cap B
\]
and note that
\[
|\tilde{S}| \leq (a\wedge r)^d
\]
where $r$ is the side of $B$.  We will estimate 
$P(A(n,m,M,S))$ for $M=2^k$, $k\geq0$.  There are two cases.  

\noindent
\textbf{Case 1:}  $n\ne m$.  

According to Lemma \ref{poisson-est},
\begin{eqnarray*}
P(A(n,m,M,S)) &=& P(N(\tilde{S},n)\ge 1)P(N(\tilde{S},m)\ge M)   \\
&\le& Ce\lambda(\tilde{S},n)\frac{\lambda(\tilde{S},m)^M e^M}{M^M}  \\
&=& Ce|\tilde{S}|^{M+1}\frac{2^{np+mMp}e^M}{M^M}    \\
&\le& Ce\cdot(a\wedge r)^{d(M+1)}\frac{2^{np+mMp}e^M}{M^M}.
\end{eqnarray*}
\noindent
\textbf{Case 2:}  $n=m$.  
\begin{eqnarray*}
P(A(n,m,M,S)) &=& P(N(\tilde{S},m)\ge M+1)   \\
&\le& C\frac{\lambda(\tilde{S},n)^{M+1} e^M}{M^M}  \\
&=& C|\tilde{S}|^{M+1}\frac{2^{np+mMp}e^M}{M^M}  \\
&\le& C(a\wedge r)^{d(M+1)}\frac{2^{np+mMp}e^M}{M^M}.
\end{eqnarray*}

Next, for some $p_0>0$, let 
\[
a_{n,m,M} = e^{-1/d}\delta^{4/(Md)} 2^{-pn/(dM)} 2^{-pm/d} M^{1/d}
 2^{-\delta(n+m+k)/(dM)}.
\]
Note that
\[
a_{n,m,M}^{dM} = \delta^4 
\left[\frac{2^{np+mMp}e^M}{M^M}\right]^{-1}
 2^{-\delta(n+m+k)}
\]
and therefore, for $S\in\mathbf{C}_{3a_{n,m,M}}(B,a)$,
\begin{eqnarray*}
\lefteqn{ (a_{n,m,2^k}\wedge r)^{-d}P(A(n,m,M,S))   }\\
&\hspace{.5in}\leq& (a_{n,m,2^k}\wedge r)^{-d}(a_{n,m,2^k}\wedge r)^{d(M+1)}
    \frac{2^{np+mMp}e^M}{M^M}  \\
&\hspace{.5in}=& (a_{n,m,2^k}\wedge r)^{dM}
    \frac{2^{np+mMp}e^M}{M^M}  \\
&\hspace{.5in}\leq& \delta^4 2^{-\delta(n+m+k)}.
\end{eqnarray*}
Let $A(B)$ be the event that there exist $n,m\geq-N$ and $M=2^k$ with $k\geq0$, 
such that there exists a particle of type $n$ in $B$, and within distance 
$a_{n,m,M}$ of this particle there are $M$ particles of type $m$ in $B$.  Note 
that the number of cubes in $\mathbf{C}_{3a}(B,a)$ is the same as the number 
of cubes in $\mathbf{C}_{a}(B,a)$.  This number is less than or equal to a 
constant times $(a\wedge r)^{-d}$.  Therefore, we have the following.  
We write $\sum_{S}$ for the sum over 
$S\in\mathbf{C}_{3a_{n,m,2^{k}}}(B,a_{n,m,2^{k}})$.  As before, 
$\tilde{S}=S\cap B$.  
\begin{eqnarray*}
\lefteqn{ P(A(B))  }\\
&\le& \sum_{k=0}^{\infty}\sum_{n,m=-N}^{\infty} 
  \sum_{S} P(A(n,m,M,S))                            \\
&=& \sum_{k=0}^{\infty}\sum_{n,m=-N}^{\infty} (a_{n,m,2^k}\wedge r)^{d}
  \sum_{S}(a_{n,m,2^k}\wedge r)^{-d} P(A(n,m,M,S))  \\
&\le& C\sum_{k=0}^{\infty}\sum_{n,m=-N}^{\infty} (a_{n,m,2^k}\wedge r)^{d}
  \sum_{S}\delta^4 2^{-\delta(n+m+k)}    \\
&=& C\sum_{k=0}^{\infty}\sum_{n,m=-N}^{\infty} (a_{n,m,2^k}\wedge r)^{d}
  \delta^4 2^{-\delta(n+m+k)}\#S    \\
&=& C \sum_{k=0}^{\infty}\sum_{n,m=-N}^{\infty} 
 \delta^4 2^{-\delta(n+m+k)}    \\
&=& C \delta^4\sum_{k=0}^{\infty}2^{-\delta k}
\left( \sum_{n=-N}^{\infty} 2^{-\delta n} \right)^2 \\
&=& C \delta^4 \frac{1}{1-2^{-\delta}}
\left( \frac{2^{\delta N}}{1-2^{-\delta}} \right)^2 \\
&\le& C \delta
\end{eqnarray*}
where $C$ depends on $N$, but not on $\delta$, as long as $\delta<1$.  Here, 
$\#S$ denotes the number of cubes $S$ that we are considering.  

Thus, we have shown:
\begin{lemma}
\label{distance}
There exists a constant $C>0$ such that the following holds.  
Fix a ball $B\subset\mathbf{R}^d$ of radius at least 1, let $\delta\in(0,1)$, 
and let
\[
a_{n,m,M} = e^{-1/d}\delta^{4/(Md)} 2^{-pn/(dM)} 2^{-pm/d} M^{1/d}
 2^{-\delta(n+m+k)/(dM)}.
\]
Let $A$ be the event that for some $n,m,M$, there exists a particle of type $n$ 
in $B$ with $M$ particles of type $m$ in $B$ within distance $a_{n,m,M}$.  Then
\[
P(A) \le C\delta.
\]
\end{lemma}

\subsection{Estimation of the solution for the heat equation.}

Now we will use Lemma \ref{distance} to bound the solution $u(t,x)$.  
Assume that event $A$ does not occur, where this event was defined in Lemma 
\ref{distance}.  

Taking $x$ as fixed, define $A_1$ to be the event that for some $m,M$, there 
exist $M$ particles of type $m$ within distance $a_{0,m,M}$ of $x$.  The same 
arguments as before easily lead to the estimate
\[
P(A_1) \leq C\delta
\]
where $C$ is the same constant as in Lemma \ref{distance}.  

For the rest of the section, fix $n$, and consider $u_n(t,x)$.  Here $n$ is 
no longer a subscript for $a_{n,m,M}$.  We will assume that neither 
$A_1$ nor $A$ occur.  Also, we write
\begin{eqnarray*}
x^{(n+1)} &=& x      \\
M &=& 2^k.  
\end{eqnarray*}
Define
\[
\mathcal{M}(n) = \left\{
(i_1,\dots,i_n)\in\mathbf{N}^n\; :\; \forall 1\leq k\leq n-1,\; i_k\ne i_{k+1}
\right\}.
\]
Let $v(t,\cdot)$ be the solution of the heat equation on $D$ with initial 
function $u_0$, and recall that $G(t,x,y)\leq G(t,x-y)$.  Setting 
$i_{n+1}=n+1$ and $m_{n+1}=0$, we get
\begin{eqnarray}
\label{7.2}
\lefteqn{ |u_n(t,x)|\mathbf{1}(A^c\cap A_1^c) }     \\
&\leq& \sum_{(i_1,\dots,i_n)\in\mathcal{M}(n)} 
\int_{\mathbf{T}_n}     \nonumber\\
&& \hspace{.5in}\prod_{j=1}^{n}\left[G(s_{j+1}-s_{j},
           x^{(i_{j+1})}-x^{(i_j)})y_{i_j} \right] |v(s_1,x^{(i_1)})|ds   \nonumber\\
&\le& \|u_0\|_\infty
\sum_{k_1,\dots,k_n=0}^{\infty}
\sum_{m_1,\dots,m_n=-N}^{\infty}  
\int_{\mathbf{T}_n}   \nonumber\\
&& \hspace{.5in}\prod_{j=1}^{n}
\left[G(s_{j+1}-s_{j},a_{m_{j+1},m_j,2^{k_j}})2^{k_j+1}2^{-m_j} \right]ds.  \nonumber
\end{eqnarray}
We assumed that $\|u_0\|_\infty\leq 1$.  It is known that the 
solution of the heat equation is bounded when the initial function $u_0(x)$ is 
bounded, and $\|v(t,\cdot)\|_\infty\leq\|u_0\|_\infty$.  When summing over all 
$(i_1,\dots,i_n)\in\mathcal{M}_n$ we sum first over all possible types 
$m_1,\dots,m_n$ of particles $x^{(i_1)},\dots,x^{(i_n)}$.  For two fixed 
consecutive particles $x^{(i_j)}$ and $x^{(i_{j+1})}$, we let $k_j$ be the 
smallest nonnegative integer $k$ such that 
\[
a_{m_{j+1},m_j,2^k}
\leq
\left|x^{(i_{j+1})}-x^{(i_j)}\right| 
<
a_{m_{j+1},m_j,2^{k+1}}
\]
where $m_{j+1}$ and $m_j$ are the types of the particles $x^{(i_{j+1})}$ and 
$x^{(i_j)}$.  Observe that such a $k_j$ always exists.  Because we are on the 
complement of both sets $A$ and $A_1$, there is no particle of type $m_j$ at a 
distance less than $a_{m_{j+1},m_j,0}$ from $x^{(i_{j+1})}$.  Thus, 
\[
\left|x^{(i_{j+1})}-x^{(i_j)}\right|\geq a_{m_{j+1},m_j,2^0}.  
\]
On the other hand, for fixed $m_{j+1}$ and $m_j$, we have
\[
\lim_{k\to\infty}a_{m_{j+1},m_j,2^k} = \infty.  
\]
Thus for large enough $k$, we have
\[
\left|x^{(i_{j+1})}-x^{(i_j)}\right|\geq a_{m_{j+1},m_j,2^k}.  
\]
Therefore,
\[
G(s_{j+1}-s_j,x^{(i_{j+1})}-x^{(i_j)}) \leq G(s_{j+1}-s_j,a_{m_{j+1},m_j,2^k}).
\]
Since $x^{(i_j)}$ and $x^{(i_{j+1})}$ are within distance 
$a_{m_{j+1},m_j,2^{k_j+1}}$ for a fixed $k_j$ and a fixed $x^{(i_{j+1})}$, there
could be at most $2^{k_j+1}-1<2^{k_j+1}$ particles $x^{(i_j)}$ of type $m_j$ 
within this distance from $x^{(i_{j+1})}$.  This is the reason why the factor 
$2^{k_j+1}$ appears in the previous inequality.  Since $x^{(i_j)}$ is of type 
$m_j$, we have $y_{i_j}\leq2^{-m_j}$ and this is why $2^{-m_j}$ appears in 
(\ref{7.2}).

Estimating the term $G$ in (\ref{7.2}), with $M=2^k$, we get 
\begin{eqnarray*}
\lefteqn{
G(r,a_{m_{j+1},m_j,2^{k}})  }\\
&=& Cr^{-d/2}\exp\left(-\frac{a^2_{m_j,m_{j+1},M}}{4r}\right)     \\
&=& Cr^{-d/2}\exp\bigg(-(4r)^{-1}
e^{-2/d}\delta^{8/(Md)} 2^{-2pm_{j+1}/(dM)} 2^{-2pm_j/d}  \\
&&\hspace{1.5in}  \cdot M^{2/d} 2^{-2\delta(m_{j+1}+m_j+k)/(Md)} \bigg)     \\
&\le& Cr^{-d/2}\exp\bigg(-(4r)^{-1}
e^{-2/d}\delta^{8/d}        \\
&& \hspace{1.5in} 2^{2[k(M-\delta)-m_{j+1}(p+\delta)-m_j(pM+\delta)]/(dM)}
 \bigg).
\end{eqnarray*}
Therefore,
\[
G(r,a_{m_{j+1},m_j,2^{k}}) 
\le Cr^{-d/2}\exp\left(-C_0r^{-1}
2^{2[k(1-\delta)-(m_{j}+m_{j+1}/M)(p+\delta)]/d}
 \right)
\]
where $C_0$ depends on $\delta$.  Thus,
\begin{eqnarray*}
\lefteqn{ u_n(t,x) 
\le \int_{\mathbf{T}_n}   
\sum_{k_1,\dots,k_n=0}^{\infty}
\sum_{m_1,\dots,m_n=-N}^{\infty}             }\\
&&  \prod_{j=1}^{n}\Bigg[C\left(s_{j+1}-s_{j}\right)^{-d/2} 
        2^{k_j-m_j}             \\
&&\cdot \exp\left(-C_0\left(s_{j+1}-s_{j}\right)^{-1}
2^{2[k_j(1-\delta)-(m_{j}+m_{j+1}/M)(p+\delta)]/d}
\right) \Bigg]ds.
\end{eqnarray*}
Rearranging the terms and using $2^{-m_1/(M+1)}\leq1$ and 
since $m_{n+1}=0$, we get
\begin{eqnarray}
\label{7.3}
\lefteqn{ u_n(t,x)
\le
\int_{\mathbf{T}_n}   
\sum_{k_1,\dots,k_n=0}^{\infty}
\sum_{m_1,\dots,m_n=-N}^{\infty}             }\\
&& \prod_{j=1}^{n}\Bigg[
C(s_{j+1}-s_j)^{-d/2} 2^{k_j-(M/(M+1))(m_{j}+m_{j+1}/M)}      \nonumber\\
&&\cdot \exp\left(-C_0\left(s_{j+1}-s_{j}\right)^{-1}
2^{2[k_j(1-\delta)-(m_{j}+m_{j+1}/M)(p+\delta)]/d}
\right) \Bigg]ds.  \nonumber
\end{eqnarray}
We wish to replace the above sums over $m_j$ by integrals.  This would involve 
replacing each variable $m_j$ by a continuous variable $w_j$ taking values in 
$[m_j,m_j+1]$.   This would decrease the exponential, but would increase the 
first power of 2 by a multiplicative factor of $2^n$.  We would have
\begin{eqnarray}
\label{7.4}
\lefteqn{ |u_n(t,x)|
\le
\int_{\mathbf{T}_n}   
\sum_{k_1,\dots,k_n=0}^{\infty}
\int_{-N}^{\infty}\dots\int_{-N}^{\infty}
}\\
&& \prod_{j=1}^{n}\Bigg[
C(s_{j+1}-s_j)^{-d/2} 2^{k_j-(M/(M+1))(w_{j}+w_{j+1}/M)}      \nonumber\\
&&\cdot \exp\left(-C_0\left(s_{j+1}-s_{j}\right)^{-1}
2^{2[k_j(1-\delta)-(w_{j}+w_{j+1}/M)(p+\delta)]/d}
\right) \Bigg]dwds  \nonumber
\end{eqnarray}
where, as with $ds$, we let $dw=dw_1\dots dw_n$.  Next, we wish to make the
change of variables
\[
v_j = w_j + \frac{w_{j+1}}{M}, \qquad j=1,\dots,n.  
\]
Now $w_{n+1}$ is fixed, and the reader can easily check that the Jacobian 
determinant for this change of variables is 1.  Since 
$M=2^{k}\geq1$ and $w_j\geq-N$, we have that $v_j\geq-2N$.  After making this 
change of variables, we can put the integral signs over $v$ inside the product 
as follows.  
\begin{eqnarray}
\label{7.5}
\lefteqn{ |u_n(t,x)|
\le
\int_{\mathbf{T}_n}   
\sum_{k_1,\dots,k_n=0}^{\infty}
\int_{-2N}^{\infty}\dots\int_{-2N}^{\infty}
}\\
&& \prod_{j=1}^{n}\Bigg[
C(s_{j+1}-s_j)^{-d/2} 2^{k_j-(M/(M+1))v_j}      \nonumber\\
&&\cdot \exp\left(-C_0\left(s_{j+1}-s_{j}\right)^{-1}
2^{2[k_j(1-\delta)-v_j(p+\delta)]/d}
\right) \Bigg]dvds  \nonumber   \\
&=& \int_{\mathbf{T}_n}
   \prod_{j=1}^{n}\Bigg[\sum_{k_j=0}^{\infty}\int_{-2N}^{\infty}
C(s_{j+1}-s_j)^{-d/2} 2^{k_j-(M/(M+1))v_j}      \nonumber\\
&&\cdot \exp\left(-C_0\left(s_{j+1}-s_{j}\right)^{-1}
2^{2[k_j(1-\delta)-v_j(p+\delta)]/d}
\right)dv_j \Bigg]ds  \nonumber   
\end{eqnarray}
where once more, $dv=dv_1\dots dv_n$.  Let $T_{j,k_j}$ denote the term inside 
the final sum in (\ref{7.5}).  Dropping the subscripts on 
$T,s,k$ and setting $r=s_{j+1}-s_j$, we can write
\begin{eqnarray}
\label{7.6}
T &=& \int_{-2N}^{\infty} Cr^{-d/2} 2^{k-(M/(M+1))v}  \\
&&  \qquad  \cdot\exp\left(-C_0 r^{-1}
          2^{2[k(1-\delta)-v(p+\delta)]/d} \right)dv.  \nonumber
\end{eqnarray}
We will split the integral into two ranges for the variable $v$.  

For some values of $v$, the exponential in 
(\ref{7.6}) is close to 1.  This will happen if the terms following 
$C_0$ in (\ref{7.6}) are less than or equal to 1.  We must solve
\[
r^{-1}2^{2[k(1-\delta)-v(p+\delta)]/d} \leq 1.
\]
Simplifying the above expression, we get
\begin{equation}
\label{7.7}
2^{-2v(p+\delta)/d} \leq r2^{-2k(1-\delta)/d}
\end{equation}
or
\begin{equation}
\label{7.8}
2^{-\frac{M}{M+1}v} \leq r^{\frac{Md}{2(M+1)(p+\delta)}}
        2^{-\frac{M(1-\delta)}{(M+1)(p+\delta)}k}.
\end{equation}
Let $\Lambda=\Lambda(k,d,p)$ be the value of $v$ for which equality is 
attained in (\ref{7.8}).  We will split $T$ into two integrals,
$T=T^{(1)}+T^{(2)}$, where
\begin{eqnarray*}
T^{(1)} &=& \int_{\Lambda}^{\infty}I(v)dv   \\
T^{(2)} &=& \int_{-2N}^{\Lambda}I(v)dv       
\end{eqnarray*}
and $I(v)$ is the integrand in (\ref{7.6}).

\textbf{Case 1:  An estimate for $T^{(1)}$.}  

By (\ref{7.8}), and using the fact that
\[
\frac{1}{2} \leq \frac{M}{M+1} \leq 1
\]
we have
\begin{eqnarray*}
T^{(1)} &=& \int_{\Lambda}^{\infty}
  Cr^{-d/2} 2^{k-(M/(M+1))v}
  \exp\left(-C_0 r^{-1} 2^{2[k(1-\delta)-v(p+\delta)]/d} \right)dv  \\
&\leq& 
  \int_{\Lambda}^{\infty}
  Cr^{-d/2} 2^{k-(M/(M+1))v}dv   \\
&\leq& Cr^{-d/2}2^{k-(M/(M+1))\Lambda}            \\
&=& Cr^{-\frac{d}{2}+ \frac{Md}{2(M+1)(p+\delta)}}
        2^{k\left(1-\frac{M(1-\delta)}{(M+1)(p+\delta)}\right)}   \\
&\leq& Cr^{-\frac{d}{2}+ \frac{d}{4(p+\delta)}}
        2^{k\left(1-\frac{M(1-\delta)}{(M+1)(p+\delta)}\right)}.
\end{eqnarray*}
since $0<r<t$.  Here, $C$ depends on $t$.  

If
\[
p<1
\]
and if $\delta$ is sufficiently small, then (remembering that $M=2^k$),
\[
\sum_{k=0}^{\infty}2^{k\left(1-\frac{M(1-\delta)}{(M+1)(p+\delta)}\right)}
=\sum_{k=0}^{\infty}2^{k\left(1-\frac{2^k(1-\delta)}{(2^k+1)(p+\delta)}\right)}
< \infty  
\]
so that 
\begin{equation}
\label{7.9}
\sum_{k_j=0}^{\infty}T^{(1)}_{j,k_j} 
  \leq C r^{-\frac{d}{2}+ \frac{d}{4(p+\delta)}}.
\end{equation}

\textbf{Case 2:  An estimate for $T^{(2)}$.}  

Now we turn to the case $v\leq\Lambda$.  Here, we assume that
\[
\Lambda \geq -2N
\]
or there will be no values of $v$ for which $-2N\leq v\leq \Lambda$.  We will 
again drop the subscripts and set $r=s_{j+1}-s_j$.  Since equality is attained 
in (\ref{7.7}) when $v=\Lambda$, we have
\begin{eqnarray*}
\lefteqn{  T^{(2)} }\\
&=& 
\int_{-2N}^{\Lambda}
  Cr^{-d/2} 2^{k-(M/(M+1))v}
  \exp\left(-C_0 r^{-1} 2^{2[k(1-\delta)-v(p+\delta)]/d} \right)dv  \\
&=& 
\int_{-2N}^{\Lambda}
Cr^{-d/2}2^{k-(M/(M+1))v}            
\exp\left(-C_0
2^{2(\Lambda-v)(p+\delta)/d} \right)dv.
\end{eqnarray*}
Changing variables to $z=\Lambda-v$, and since $M/(M+1)<1$, we see that
\begin{eqnarray*}
T^{(2)} &=&  \int_{0}^{\Lambda+2N}
Cr^{-d/2}2^{k-(M/(M+1))(\Lambda-z)}            
\exp\left(-C_0
2^{2z(p+\delta)/d} \right)dz    \\
&\leq&
 Cr^{-d/2}2^{k-(M/(M+1))\Lambda}            
\int_{0}^{\infty} 2^{z} 
\exp\left(-C_0 2^{2z(p+\delta)/d} \right)dz      \\
&\leq& 
Cr^{-d/2}2^{k-(M/(M+1))\Lambda}            \\
&\leq& Cr^{-\frac{d}{2}+\frac{d}{4(p+\delta)}}
2^{k\left(1-\frac{M(1-\delta)}{(M+1)(p+\delta)}\right)}.
\end{eqnarray*}
Once again, if
\[
p < 1
\]
and $\delta$ is small, then
\[
\sum_{k=0}^{\infty}2^{k\left(1-\frac{M(1-\delta)}{(M+1)(p+\delta)}\right)}
= \sum_{k=0}^{\infty}2^{k\left(1-\frac{2^k(1-\delta)}{(2^k+1)(p+\delta)}\right)}
< \infty
\]
so that 
\begin{equation}
\label{7.10}
T^{(2)} \leq Cr^{-\frac{d}{2}+\frac{d}{4(p+\delta)}}.
\end{equation}

Combining (\ref{7.9}) and (\ref{7.10}), if 
\[
p < 1
\]
and $\delta$ is sufficiently small, then
\[
T = T^{(1)} + T^{(2)} \leq Cr^{-\frac{d}{2}+\frac{d}{4(p+\delta)}}
\]
and
\[
u_n(t,x) \le 
C\int_{\mathbf{T}_n}
\prod_{i=1}^{n} \left(s_{i+1}-s_i\right)^{-\frac{d}{2}+\frac{d}
{4(p+\delta)}}
\;ds.
\]

Now, a calculation of 
Hu (\cite{hu01}, proof of Theorem 4.1, equation (4.10)) 
gives, for $\beta>0$,
\begin{equation}
\label{7.11}
\int_{\mathbf{T}_n}\prod_{k=1}^{n} (s_{k+1}-s_k)^{-\beta}ds
= \frac{t^{n(1-\beta)}\Gamma^n(1-\beta)}
{\Gamma\left(1+n\left(1-\beta\right)\right)}.
\end{equation}
We give a quick derivation of Hu's estimate in the appendix, since he merely 
cites the result.  

In (\ref{7.11}), substitute
\[
\beta = \frac{d}{2}-\frac{d}{4(p+\delta)}.
\]
If $\beta<1$, then $\Gamma(1-\beta)$ is defined.  With $\mathcal{C}_n$ 
denoting a term of order $C^n$ for some $C>0$, we get
\[
u_n(t,x) \le 
\frac{t^{n(1-\beta)}\Gamma^n(1-\beta)}
{\Gamma\left(1+n\left(1-\beta\right)\right)}
\leq   \mathcal{C}_n t^{n(1-\beta)} n^{-n\left(1-\beta\right)}
\]
since only the denominator has order greater than $\mathcal{C}_n$.  
Thus, if $\beta<1$, and if $A$ and $A_1$ do not occur,
\[
u(t,x) = \sum_{n=0}^{\infty} u_n(t,x) < \infty.  
\]
Solving for $\beta<1$, we obtain the requirement that
\[
\frac{d}{2}-\frac{d}{4(p+\delta)} < 1.
\]
For small enough values of $\delta$, this is satisfied if either $d\leq4$ or
\[
p < \frac{1}{2} + \frac{1}{d-2}, \qquad d\geq5.
\]
This completes the proof of Theorem \ref{t1}, for $L_K$.  As mentioned 
earlier, since $P(\mathcal{A}_K)=P(L=L_K)\to 1$ as $K\to\infty$, it follows 
that Theorem \ref{t1} holds almost surely.  

\section{Proof of the existence part of \\ \mbox{}\hspace{0.8cm} Theorem \ref{t2} ($p\geq1$)}
\label{s9}
\setcounter{equation}{0}

For $p\geq1$, we also wish to replace $L$ by $L_K$.  But in this case our 
stochastic integrals will no longer be martingales unless we delete the 
compensator for $L-L_K$.  We define this compensator as follows.  
\[
\mathcal{L}(dtdx) = dtdx\int_{K}^{\infty}c_py^{-(p+1)}dy
= c'_pK^{-p}dtdx.  
\]
Then, on the event $\mathcal{A}_K$, we can replace (\ref{1.1}) by
\begin{eqnarray*}
\frac{\partial u}{\partial t} &=& \Delta u + c'_pK^{-p}u
        + u\diamond\dot{L}_K(x)  \\
u(0,x) &=& u_0(x).  
\end{eqnarray*}
If we let 
\[
w(t,x) = \exp\left(-c'_pK^{-p}t\right)u(t,x)
\]
then $w(t,x)$ solves
\begin{eqnarray}
\label{8.1}
\frac{\partial u}{\partial t} &=& \Delta u + u\diamond\dot{L}_K(x)  \\
u(0,x) &=& u_0(x).  \nonumber
\end{eqnarray}
Therefore, any conclusions which hold almost surely for $w(t,x)$, for all $K$, 
must hold for solutions to (\ref{1.1}).

For ease of notation, we write $u(t,x)$ instead of $w(t,x)$ and $L$ instead of 
$L_K$.  With this notation, we can prove the following lemma.
\begin{lemma}
\label{poissonian_Wiener_Ito} For any $q \in (p, 2]$
there exists a positive constant $C_{q}$ such that
for any measurable symmetric function $f : D^{n} \to {\bf R}$
we have: 
\begin{eqnarray}
\label{8.2}
\lefteqn{ E\left[\left|n!\int_{\widetilde{D}_{n}}f(x^{(1)}, \dots, x^{(n)})
L(dx^{(1)}) \dots L(dx^{(n)})\right|^{q}\right] }\\
& \leq &C_{q}^{n}(n!)^{q - 1}
\int_{D^{n}}|f(x^{(1)}, \dots, x^{(n)})|^{q}dx^{(1)} \dots dx^{(n)}. \nonumber
\end{eqnarray}
\end{lemma}
Note that the above lemma gives sufficient conditions for the existence of the 
multiple stochastic integral $I_n(f_n)$.  

\textbf{Proof.} 
We would like to use the Burkholder-Davis-Gundy inequality, and for this reason 
we need to expand on our formulation of the multiple stochastic integral 
given in Section (\ref{precise}), so that we have a martingale. 
We start with the $\sigma$-algebra.  Recall that we have chosen
the first coordinate of $x$ as our time variable.  In particular, we can
define a filtration.  Let $m(D),M(D)$ be the minimum and maximum 
values of $x_1$, the first coordinate of $x$, for $x\in D$.  For 
$z\in[m(D),M(D)]$, we define the $\sigma$-algebra $\mathcal{H}_z$ as follows.  
\[
\mathcal{H}_z = \sigma \bigg\{L(A):  A\subset \{x\in D: x_1\leq z\}\bigg\}.
\]
Next, let
\[
R(z,D) = \{x\in D : x_1 \leq z\}.
\]
Let $h(x)$ be a predictable function with respect to the $\mathcal{H}_z$ filtration.  Let
\[
Y_z = \int_{R(z,D)}h(x)L(dx).
\]
Then $(Y_z,\mathcal{H}_z)$ is a martingale, for appropriate conditions on $h$, 
and its quadratic variation is
\[
\langle Y \rangle_z = \sum_{x_1^{(i)}\leq z} h^2(x^{(i)}) 
y_i^2.  
\]
For ease of notation, we drop the subscript $z$ when $z=M(D)$, so that
\begin{eqnarray*}
Y &=& Y_{M(D)},  \\
\langle Y \rangle &=& \langle Y \rangle_{M(D)}.
\end{eqnarray*}
The Burkholder-Davis-Gundy inequality now implies that for $q>1$, there is a 
constant $C_q$ depending only on $q$, such that
\[
E\left[\left|Y\right|^q\right]
\leq C_q E\left[\langle Y\rangle^{q/2}\right].
\]
Recall that for $r\leq1$ and a sequence of nonnegative numbers $a_i$, the 
following elementary inequality holds.  
\[
\left(\sum_{i=1}^{\infty}a_i\right)^r 
\leq \sum_{i=1}^{\infty} a_i^r.
\]
Setting $r=q/2$, and using our formula for $\langle Y\rangle$, we get
\begin{eqnarray}
\label{8.3}
E\left[\left|Y\right|^q\right]
&\leq& C E\left[\left(\sum_{i=1}^\infty
h^2(x^{(i)})y_i^2\right)^{q/2}\right]  \\
&\leq& C E\left[\sum_{i=1}^\infty
\left|h^2(x^{(i)})\right|^{q/2}y_i^{q}\right]  \nonumber\\
&=& C \int_{D}
E\left|h(x)\right|^q dx  \nonumber
\end{eqnarray}
provided $q>p$. The constant $C$ depends on $p$, $q$, and $K$. 
Applying (\ref{8.3}) repeatedly, we get
\begin{eqnarray}
\label{8.4}
\lefteqn{ E\left[\left|\int_{\widetilde{D}_n}f(x^{(1)}, \dots, x^{(n)})
L(dx^{(1)}) \dots L(dx^{(n)})\right|^q\right] }\\
& \leq & C^{n} \int_{\widetilde{D}_n}
\left|f(x^{(1)}, \dots, x^{(n)})\right|^{q}
dx^{(1)}\dots dx^{(n)}. \nonumber
\end{eqnarray}
Note, that the repeated application of inequality (\ref{8.3}) was possible due 
to the order $x_{1}^{(1)} < x_{1}^{(2)} < \dots < x_{1}^{(n)}$ that we imposed on 
$\widetilde{D}_n$. To see this, let us remember that the function $h$ involved 
in (\ref{8.3}) needed to be predictable. Thus,
\begin{eqnarray*}
\lefteqn{ E\left[\left|\int_{\widetilde{D}_n}f(x^{(1)}, \dots, x^{(n)})
L(dx^{(1)}) \dots L(dx^{(n)})\right|^q\right] }\\
& = & E\left[\left|\int_{D}h(x^{(n)})L(dx^{(n)})\right|^q\right],
\end{eqnarray*}
where 
\[
h(x^{(n)}) = \int_{\tilde{R}_{n-1}(x_{1}^{(n)})}f(x^{(1)}, \dots, x^{(n - 1)}, x^{(n)})
L(dx^{(1)}) \dots L(dx^{(n - 1)})
\]
and 
\begin{equation}
\label{8.5}
\tilde{R}_{n-1}(y)=\left\{(x^{(1)},\dots,x^{(n-1)})\in\widetilde{D}_{n-1}:
x_1^{(n-1)}<y\right\}
\end{equation}
where $y\in\mathbf{R}$.  
The function $h(x^{(n)})$ is predictable since $x_{1}^{(1)} < x_{1}^{(n)}$, $\dots$, 
$x_{1}^{(n - 1)} < x_{1}^{(n)}$ on $\tilde{R}(x_{1}^{(n)}, D)$. 
So, we can apply the inequality (\ref{8.3}) and obtain 
\begin{eqnarray*}
\lefteqn{ E\left[\left|\int_{\widetilde{D}_n}f(x^{(1)}, \dots, x^{(n)})
L(dx^{(1)}) \dots L(dx^{(n)})\right|^q\right] }\\
& \leq & 
C\int_{D}E\left[\left|h(x^{(n)})\right|^q\right]dx^{(n)},
\hspace{1in}
\end{eqnarray*}
and so on.  Using (\ref{8.3}), and the fact that $f$ is symmetric, we get:
\begin{eqnarray*}
\lefteqn{ E\left[\left|n!\int_{\widetilde{D}_{n}}f(x^{(1)}, \dots, x^{(n)})
L(dx^{(1)}) \dots L(dx^{(n)})\right|^q\right] }\\
& = & (n!)^{q}E\left[\left|\int_{\widetilde{D}_n}f(x^{(1)}, \dots, x^{(n)})
L(dx^{(1)}) \dots L(dx^{(n)})\right|^q\right] \\
& \leq & (n!)^{q}C^{n} \int_{\widetilde{D}_n}
\left|f(x^{(1)}, \dots, x^{(n)})\right|^{q}
dx^{(1)}\dots dx^{(n)}\\
& = & (n!)^{q}C^{n} \frac{1}{n!}\int_{D^{n}}
\left|f(x^{(1)}, \dots, x^{(n)})\right|^{q}
dx^{(1)}\dots dx^{(n)}\\
& = & C^{n}(n!)^{q - 1}\int_{D^{n}}
\left|f(x^{(1)}, \dots, x^{(n)})\right|^{q}
dx^{(1)}\dots dx^{(n)}.
\end{eqnarray*}
\qed\\

Before proving Theorem \ref{t2}, we present two simple lemmas.  
\begin{lemma} Let $n$ be a positive integer and $q$ a real number, such that $q \geq 1$.
Let $f : D^{n} \to {\bf R}$ be a function in $\mathbf{L}^{q}(D^{n}, dx)$, 
where $dx$ is Lebesgue measure. If $\tilde{f}$ is the symmetrization of $f$, 
then $\tilde{f} \in \mathbf{L}^{q}(D^{n}, dx)$
and 
\begin{eqnarray}
\label{8.6}
\parallel \tilde{f} \parallel_{q} & \leq & \parallel f \parallel_{q}
\end{eqnarray}
where 
$\parallel \cdot \parallel_{q}$ denotes the norm of the space 
$\mathbf{L}^{q}(D^{n}, dx)$.
\end{lemma}

\textbf{Proof.} For any permutation $\sigma$ of the set $\{1, 2, \dots, n\}$, we can consider 
the function $f_{\sigma} : D^{n} \to {\bf R}$, defined by
\[
f_{\sigma}(x^{(1)}, \dots, x^{(n)}) = f(x^{(\sigma(1))}, \dots, x^{(\sigma(n))}).
\]
Making the change of variable $y^{(1)} = x^{(\sigma(1))}$, $\dots$, 
$y^{(n)} = x^{(\sigma(n))}$,
we can see that
\begin{eqnarray*}
\lefteqn{ \int_{D^{n}}|f_{\sigma}(x^{(1)}, 
    \dots ,x^{(n)})|^q dx^{(1)} \dots dx^{(n)}     }\\
& = & \int_{D^{n}}|f(x^{(1)}, \dots ,x^{(n)})|^q dx^{(1)} \dots dx^{(n)}.
\end{eqnarray*}
Thus for all permutations $\sigma$, we have 
\begin{eqnarray*}
\parallel f_{\sigma} \parallel_{q} & = & \parallel f \parallel_{q}.
\end{eqnarray*}
Since the symmetrization $\tilde{f}$ of $f$ is defined as
\[
\tilde{f} = \frac{1}{n!}\sum_{\sigma}{f_{\sigma}},
\]
we can apply the triangle inequality to get
\begin{eqnarray*}
\parallel f_{\sigma} \parallel_{q} & = & 
\parallel \frac{1}{n!}\sum_{\sigma} f_{\sigma} \parallel_{q}\\
& \leq & \frac{1}{n!}\sum_{\sigma}\parallel f_{\sigma} \parallel_{q}\\
& = & \frac{1}{n!}\sum_{\sigma}\parallel f \parallel_{q}\\
& = & \parallel f \parallel_{q}.
\end{eqnarray*}

\qed

In the proof of Lemma \ref{poisson-est} we used the inequality 
$e^{n} \geq n^{n}/n!$, for all positive integers $n$. We would like to prove 
this inequality for all positive real numbers $x$ and therefore we have to 
replace $n!$ by $\Gamma(x + 1)$. Thus we will prove the following lemma.

\begin{lemma}
\label{gamma}
For any positive real number $x$, we have
\begin{eqnarray}
\label{8.7}
\Gamma(x + 1) & \geq & x^{x}e^{-x}.
\end{eqnarray}
\end{lemma}

\textbf{Proof.} If $x > 0$, then for all $t \geq x$, $t^{x} \geq x^{x}$. Thus we have: 

\begin{eqnarray*}
\Gamma(x + 1) & = & \int_{0}^{\infty} t^{x}e^{-t}dt\\
& \geq & \int_{x}^{\infty} t^{x}e^{-t}dt\\
& \geq & \int_{x}^{\infty} x^{x}e^{-t}dt\\
& = & x^{x}\int_{x}^{\infty} e^{-t}dt\\
& = & x^{x}e^{-x}.
\end{eqnarray*}

\qed\\

We can prove now Theorem \ref{t2}.
According to inequalities (\ref{8.2}) and (\ref{8.6}), we have:
\begin{eqnarray}
\label{8.8}
\lefteqn{   E\left[\left|u_n(t,x)\right|^q\right]   }\\
&\leq& C^{n}(n!)^{q - 1}\int_{D^{n}} \bigg[\;
        {\rm sym}_{x^{(1)},\dots, x^{(n)}}|u_0(x^{(1)})|   \nonumber\\
& \ &  \hspace{1in}        \cdot\int_{\mathbf{T}_n}
   \prod_{k=1}^{n} G(s_{k+1}-s_k,x^{(k+1)}-x^{(k)})ds\bigg]^q dx   \nonumber\\
&\leq& aC^{n}(n!)^{q - 1}\int_{D^{n}} \left[\;
        \int_{\mathbf{T}_n}
   \prod_{k=1}^{n} G(s_{k+1}-s_k,x^{(k+1)}-x^{(k)})ds\right]^q dx,     \nonumber
\end{eqnarray}
where $a = \parallel u_{0} \parallel_{\infty}$ and 
\[
\mathbf{T}_{n} = \{(s_{1}, \dots, s_{n}) \in \mathbf{R}^{n}: 
 0 < s_{1} < s_{2} < \dots <s_{n} <t\}.
\]
We may assume that $\parallel u_{0} \parallel_{\infty} = 1$.
Let $V(\mathbf{T}_n)=t^n/n!$ be the volume of $\mathbf{T}_n$.  Using Jensen's 
inequality, we obtain
\begin{eqnarray}
\label{8.9}
\lefteqn{  \left[\int_{\mathbf{T}_n}
     \prod_{k=1}^{n} G(s_{k+1}-s_k,x^{(k+1)}-x^{(k)})ds\right]^q    }   \\
&=& V(\mathbf{T}_n)^q\left[V(\mathbf{T}_n)^{-1}\int_{\mathbf{T}_n}
     \prod_{k=1}^{n} G(s_{k+1}-s_k,x^{(k+1)}-x^{(k)})ds\right]^q       \nonumber\\
&\leq& V(\mathbf{T}_n)^{q}V(\mathbf{T}_n)^{-1}\int_{\mathbf{T}_n}
     \prod_{k=1}^{n} G^q(s_{k+1}-s_k,x^{(k+1)}-x^{(k)}) ds     \nonumber\\
&=& \left(\frac{t^{n}}{n!}\right)^{q-1}   \int_{\mathbf{T}_n}
     \prod_{k=1}^{n} G^q(s_{k+1}-s_k,x^{(k+1)}-x^{(k)}) ds.       \nonumber
\end{eqnarray}
Recall that
\begin{eqnarray}
\label{8.10}
\int_{\mathbf{R}^d} G^q(r,x) dx
&=& Cr^{-(q-1)d/2}
  \int_{\mathbf{R}^d} r^{-d/2}\exp\left(-\frac{qx^2}{4r}\right) dx   \nonumber\\
&=& Cr^{-(q-1)d/2}.
\end{eqnarray}
Combining (\ref{8.8}) and (\ref{8.9}), we see that $(n!)^{q - 1}$ is canceled out 
by $(1/n!)^{q - 1}$, and by using Fubini's theorem and (\ref{8.10}) repeatedly (i.e. 
integrated first over $x^{(1)}$, then $x^{(2)}$, and so on), we obtain:
\begin{eqnarray}
\label{8.11}
\lefteqn{   E\left[\left|u_n(t,x)\right|^q\right]   }\\
& \leq & \left(Ct^{q - 1}\right)^{n}\int_{\mathbf{T}_{n}}\left[\int_{D^{n}}
\prod_{k=1}^{n} G^q(s_{k+1}-s_k,x^{(k+1)}-x^{(k)}) dx\right]ds \nonumber\\
& \leq & \left(Ct^{q - 1}\right)^{n}\int_{\mathbf{T}_{n}}\left[\int_{\mathbf{R}^{n}}
\prod_{k=1}^{n} G^q(s_{k+1}-s_k,x^{(k+1)}-x^{(k)}) dx\right]ds \nonumber\\
& = & \left(Ct^{q - 1}\right)^{n}\int_{\mathbf{T}_{n}}
\prod_{k=1}^{n}(s_{k + 1} - s_{k})^{-(q - 1)d/2}ds. \nonumber
\end{eqnarray}
Now, an estimate of Hu~(\ref{7.11}) gives
\[
\int_{\mathbf{T}_n}\prod_{k=1}^{n} (s_{k+1}-s_k)^{-(q-1)d/2}ds
= \frac{t^{n(1-(d/2)(q-1))}\Gamma^n(1-\frac{d}{2}(q-1))}
{\Gamma\left(1+n\left(1-\frac{d}{2}(q-1)\right)\right)}.
\]

If $d<2/(q-1)$, then $\Gamma(1-(d/2)(q-1))$ is defined.  
For $d < 2/(q - 1)$, let $\alpha = 1- (d/2)(q - 1) > 0$. 
Then from (\ref{8.11}) we conclude that
\begin{equation}
\label{8.12}
E\left[|u_n(t, x)|^q\right]  \leq  \frac{\rho^{n}}{\Gamma(1 + n\alpha)},
\end{equation}
for some constant $\rho$ depending on $t$, $q$, $p$, and $K$.
Using the well-known inequality $(b + c)^{q} \leq 2^{q - 1}(b^{q} + c^{q})$, for all $q \geq 1$ and 
$b$ and $c$ non-negative numbers, an easy induction argument shows 
that for $a_n\geq0$, and $c=2^{q-1}$,
\[
\left(\sum_{n=0}^{\infty}a_n\right)^q \leq \sum_{n=0}^{\infty}c^{n+1}a_n^{q}.
\]
Applying this inequality to $u(t,x)$, we find
\begin{eqnarray}
\label{8.13}
E\left[|u(t, x)|^q\right] & = & E\left[\left(\sum_{n = 0}^{\infty}|u_n(t, x)|\right)^q\right]\\
& \leq & E\left[\sum_{n = 0}^{\infty}c^{n + 1}|u_n(t, x)|^{q}\right] \nonumber\\
& = & \sum_{n = 0}^{\infty}c^{n + 1}E\left[|u_n(t, x)|^{q}\right] \nonumber\\
& \leq & c\sum_{n = 0}^{\infty}\frac{1}{\Gamma(1 + n\alpha)}(c\rho)^{n}. \nonumber
\end{eqnarray}
Using now (\ref{8.7}) we obtain:
\begin{eqnarray}
\label{8.14}
E\left[|u(t, x)|^q\right] & \leq & c\sum_{n = 0}^{\infty}\frac{e^{n\alpha}}{(n\alpha)^{n\alpha}}
(c\rho)^{n}.
\end{eqnarray}
Thus, we can see from (\ref{8.14}) that if 
$d < 2/(q - 1)$ (or equivalently $q < 1 + (2/d)$) and $p < q$, then 
$E\left[|u(t, x)|^q\right] < \infty$. (The condition $p < q$ was used in 
the Lemma \ref{poissonian_Wiener_Ito}.) For such a $q$ to exist we need that $p < 1 + 2/d$. \\ 
\par Let $\Pi$ be the random Poisson set $\{(x^{(i)}, y_{i})\}_{i \geq 1}$. For any positive 
integer $k$ we define the event:
\[
A_{k} = \{\omega \in \Omega ~:~ \forall~ (x, y) \in \Pi(\omega), y \leq k\}.
\]
That means, $A_{k}$ is the event that all atoms have masses less than or equal to $k$.\\

\par Thus if $p < 1 + (2/d)$, then for all $q \in (p, 1 + (2/d))$ and for all $k \geq 1$,
\[
E[|u(t, x)|1_{A_{k}}^{q}] < \infty.
\]
Since, the expectation is the integration with respect to a probability (finite) measure, 
it follows, that for any $q'$ and $q$, such that $1 \leq q' \leq p < q < 1 + (2/d)$, we have:
\[
\left(E[|u(t, x)1_{A_{k}}|^{q'}]\right)^{1/q'} \leq \left(E[|u(t, x)1_{A_{k}}|^{q}]\right)^{1/q} 
< \infty.
\]
Thus the solution $u(t, x)$ exists, for all $t \geq 0$ and all $x \in D$, and belongs to the space 
$\mathbf{L}^{q}_{loc}(\Omega)$, for all $1 \leq q < 1 + (2/d)$.

This completes the proof of Theorem \ref{t2}, if $u(t,x)$ means $w(t,x)$ and 
$L$ means $L_K$.  As mentioned at the beginning of Section \ref{s9}, this proof
carries over to $u(t,x)$ and $L$, since $P(\mathcal{A}_K)=P(L=L_K)\to1$ as 
$K\to\infty$.

\section{Uniqueness}
\setcounter{equation}{0}

\subsection{The Main Lemma}

We consider both cases $p<1$ and $p\geq1$.  Suppose $u(t,x)$ is a solution to
(\ref{1.1}) in case $p\geq1$ or (\ref{1.2}) in case $p<1$.  Then
we have
\begin{equation}
\label{9.1}
u(t,x) = \sum_{n=0}^{\infty}g_n(t,x)
\end{equation}
where $g_n\in \mathcal{J}_n$ in case $p<1$ or $g_n\in \mathcal{I}_n$ in case
$p\geq1$.

\begin{lemma}
\label{uniquely-determined}
Let $u(t,x)$ be as in (\ref{9.1}).  If the $g_n$ are uniquely
determined from $u$, then uniqueness holds for (\ref{1.1}) or 
(\ref{1.2}).
\end{lemma}

\textbf{Proof.}
Both cases have the same proof, so we will only deal with the case $p\geq1$.  
Suppose there are two solutions $u_i(t,x)$ for $i=1,2$, where
\[
u_i(t,x) = \sum_{n=0}^{\infty}g^{(i)}_n(t,x)
\]
and $g^{(i)}_n(t,x)\in\mathcal{I}_n$.  Let
\[
Y(t,x) = u_1(t,x) - u_2(t,x) = \sum_{n=0}^{\infty}g_n(t,x)
\]
where $g_n(t,x)=g^{(1)}_n(t,x)-g^{(2)}_n(t,x)$.  Then $Y(t,x)$ is also a 
solution to (\ref{1.1}) with initial condition $Y(0,x)=0$, so by 
(\ref{6.1}), we have
\begin{eqnarray*}
\sum_{n=0}^{\infty}g_n(t,x) &=& Y(t,x) \\
&=& \int_{D}\left(\int_{0}^{t}G(t-s,x,y)Y(s,y)ds\right)\diamond L(dy) \\
&=& \sum_{n=0}^{\infty}\int_{D}\left(\int_{0}^{t}G(t-s,x,y)g_n(s,y)ds\right)\diamond L(dy) \\
&=& \sum_{n=1}^{\infty}\tilde{g}_n(t,x)  \\
&&\qquad \forall  x\in \mathbf{R}^d, t\in (0,T].  
\end{eqnarray*}
where 
\begin{eqnarray}
\label{9.2}
\tilde{g}_n(t,x)
&=& \int_{D}\left(\int_{0}^{t}G(t-s,x,y)g_{n-1}(s,y)ds\right)\diamond L(dy) \\
&\in& \mathcal{I}_{n}.  \nonumber
\end{eqnarray}
Since we are assuming that $Y(t,x)$ uniquely determines both $g_n(t,x)$ and 
$\tilde{g}_n(t,x)$, it follows that $g_0(t,x)=0$ and 
\begin{equation}
\label{9.3}
g_{n}(t,x)=\tilde{g}_n(t,x)
\end{equation}
for any positive value of $n$.  Thus, (\ref{9.2}) implies that 
$\tilde{g}_1(t,x)=0$, and hence by (\ref{9.3}) $g_1(t,x)=0$.
An induction argument now yields $g_n(t,x)=0$ for all values of $n$.
\qed

\subsection{Proof of Uniqueness for $p<1$}

For simplicity in this proof we will assume from the beginning that 
 $\dot{L}$ is a stable noise (not $\dot{L}_K$) without truncation of large 
atoms. 
  
We begin by ordering the atoms of $\dot{L}$ according to size, so that
$y_1>y_2>\cdots$. 

 For $z\in(1,2)$, let $R_{i}(z)\dot{L} $ be the transformation of 
the noise $\dot{L}$ which replaces $y_i$ with $zy_i$.  

\begin{lemma}
\label{lemma-abs-cont}
Fix $z\in(1,2),\; i\geq 1$.  Let 
\[
\dot{M}=
R_i(z)\dot{L}
\]
Also, let $P_M,P_L$ be the probability measures on the canonical probability
space of atomic measures, generated by $\dot{M}$ and $\dot{L}$, respectively.
Then
\[
P_M << P_L.
\]
\end{lemma}

\textbf{Proof.}
For any $k\geq 1, n\geq 0$ define an event
\begin{eqnarray*}
 A_{k,n}&\equiv &
\left\{ \mbox{there are exactly $n$ atoms of $\dot{L}$}
\right.
\\
&&\left. \mbox{with masses greater than or equal to  $1/k$}\right\}
\end{eqnarray*}
Note that for any $k\geq 1$
\begin{eqnarray*}
\bigcup_{n=0}^{\infty} A_{k,n}&=&\Omega\\
A_{k,l}\bigcap A_{k,m} &=& \emptyset,\;\forall m\not=l.
\end{eqnarray*}
It is obvious that
\begin{eqnarray}
\label{9.4}
\lim_{k\rightarrow \infty} P_L\left( \bigcup_{l=0}^{i} A_{k,l}\right)=0.
\end{eqnarray}
Fix $\delta>0$ arbitrary small. By~(\ref{9.4}) we can fix $k$ sufficiently large such that
\begin{eqnarray}
\label{9.5}
P_L\left( \bigcup_{l=0}^{i} A_{k,l}\right)\leq \delta.
\end{eqnarray}
Define
\begin{eqnarray}
\label{9.6}
 \bar{A}_{k,i}\equiv\left(\bigcup_{l=0}^{i} 
 A_{k,l}\right)^{c}=\bigcup_{l=i+1}^{\infty} A_{k,l}.
\end{eqnarray}
Now  will show that
\begin{eqnarray}
\label{9.7}
P_{M}<< P_{L}\;\;\;\mbox{on}\; \bar{A}_{k,i}\,.
\end{eqnarray}

Let us decompose $\dot{L}$ as follows:
\[
\dot{L} = \dot{L}^{1,k} + \dot{L}^{2,k}
\]
where $\dot{L}_1$ includes atoms of size $y\geq 1/k$ and $\dot{L}_2$
includes atoms of size $y<1/k$.  Note that $\dot{L}^{1,k},\dot{L}^{2,k}$ are
independent, and
\begin{eqnarray*}
R_i(z)\dot{L}=R_i(z)\dot{L}^{1,k}+\dot{L}^{2,k}\;\;\;\;\mbox{on}\; \bar{A}_{k,i}\,.
\end{eqnarray*}
Hence, if we define
\[
\dot{M}_1=
R_i(z)\dot{L}^{1,k}
\]
then, in  order to get (\ref{9.7}), it suffices to show that
\begin{equation}
\label{9.8}
P_{M_1} << P_{L^{1,k}}\;\;\mbox{on}\;A_{k,l}\,,\;\;\forall l\geq i+1
\end{equation}
where $P_{M_1},P_{L^{1,k}}$ are the measures induced by $\dot{M}_1,\dot{L}^{1,k}$
respectively. Fix an arbitrary $l\geq i+1$.
On $A_{k,l}\,$, $\dot{M}_1$ and $\dot{L}^{1,k}$
have $l$ atoms at the same locations. Therefore, to get the desired absolute
continuity we just need to show
the absolute continuity of the laws of the corresponding atom sizes.
To be more precise, from now on we assume that $A_{k,l}$ occurs. 
Let $\hat{Y}^L_1\,,\ldots,\hat{Y}^L_l$ (resp.
$\hat{Y}^M_1\,,\ldots,\hat{Y}^M_l$)
be the unordered sizes of the atoms of $L^{1,k}$ (resp. $M_1$).
$(\hat{Y}^L_1\,,\ldots,\hat{Y}^L_l)$ and
$(\hat{Y}^M_1\,,\ldots,\hat{Y}^M_l)$ are continuous $l$-dimensional random
variables whose laws are supported on $(1/k,\infty)^{l}$. Moreover,
the probability density function of $(\hat{Y}^L_1\,,\ldots,\hat{Y}^L_l)$ is
positive at any point of $(1/k,\infty)^{l}$. This immediately implies that
conditioned on $A_{k,l}$ the probability law of
$(\hat{Y}^M_1\,,\ldots,\hat{Y}^M_l)$ is absolutely continuous with respect
to the law of $(\hat{Y}^L_1\,,\ldots,\hat{Y}^L_l)$, and since $l\geq i+1$
was arbitrary, (\ref{9.8}) and (\ref{9.7})
follow.

 The end of the proof of the lemma is trivial. Let $B$ be any event such that
 \[ P_{L}(B)=0.\]
Then
\begin{eqnarray}
\nonumber
P_{M}(B)&=& \sum_{l=0}^{i}P_{M}(B\bigcap A_{k,l}) +\sum_{l=i+1}^{\infty} P_{M}(B\bigcap A_{k,l})\\
\nonumber
&=&\sum_{l=0}^{i}P_{M}(B\bigcap A_{k,l})\\
\label{9.9}
&\leq&P_{M}(\bigcup_{l=0}^{i} A_{k,l}),
\end{eqnarray}
where the second equality follows from (\ref{9.7}). Now it is easy to check 
the following equality of events:
 \begin{eqnarray*}
\lefteqn{\left\{\mbox{there are no more than $i$ atoms of $\dot{L}$
  with masses $y\geq 1/k$}\right\}}\\
  &=&\left\{\mbox{there are no more than $i$ atoms of $\dot{M}$
  with masses $y\geq 1/k$}
\right\}
\end{eqnarray*}
Therefore, we get
\[P_{M}\left(\bigcup_{l=0}^{i} A_{k,l}\right)
  =P_{L}\left(\bigcup_{l=0}^{i} A_{k,l}\right).\]
Thus, (\ref{9.5}) and (\ref{9.9}) imply that
\[ P_M(B)\leq \delta.\]
Since $\delta>0$ was arbitrary, we get $P_M(B)=0$ and the result follows.
\qed

\mbox{}\\

By Lemma~\ref{uniquely-determined} the following lemma gives the required uniqueness for the
 case $p<1$.

\begin{lemma}
\label{unique-expansion}
Let $p<1$.
In (\ref{9.1}), the functions $g_n$ are uniquely determined.
\end{lemma}

\textbf{Proof}
We argue by contradiction.  Suppose Lemma \ref{unique-expansion} were false.
Then by subtracting the two series, we form a series
\[
{\cal G}\equiv\sum_{n=0}^{\infty}g_n = 0
\]
where $g_n\in\mathcal{J}_n$, and at least one of the $g_n$ is not identically 0.
In fact, we may assume that
\[
g_n = \sum_{(i_1,\dots,i_n)\in J_n}y_{i_1}\dots y_{i_n}h_n(x_{i_1},\dots,x_{i_n})
\]
with
\begin{equation}
\label{9.10}
\sum_{n=0}^{\infty}\sum_{(i_1,\dots,i_n)\in J_n}
  y_{i_1}\dots y_{i_n}|h_n(x_{i_1},\dots,x_{i_n})|
<\infty
\end{equation}
with probability 1.
Actually, each solution is a difference of two series for which the functions 
$h_n$ are nonnegative.

Recall that $\{y_1>y_2>\ldots\}$ are the ordered sizes of the atoms of our 
noise $\dot{L}$.  Fix an arbitrary $n\geq 1$.
Our goal is to show that, almost surely,
\[
y_1\dots y_n h_n(x_1,\dots,x_n) = 0
\]
from which it follows that $h_n=0$ almost everywhere.
By Lemma \ref{lemma-abs-cont}, for
$z_i\in (1,2)$, the probability induced by
$R_i(z_i)\dot{L}$ is absolutely continuous with respect to the
probability induced by $\dot{L}$.
 Furthermore, $R_i(z_i)$ extends to an
operator on $\mathcal{J}$, as follows. We set \[R_i(z_i)\sum_{n=0}^{\infty}g_n\equiv\sum_{n=0}^{\infty}R_i(z_i)g_n,\]
 where
we define $R_i(z_i)g_n$ to be a multiple stochastic integral with
 $\dot{L}$ replaced by $R_i(z_i)\dot{L}$.
  Therefore, if ${\cal G}=\sum_{n=0}^{\infty}g_n=0$
almost surely, then
$R_i(z_i)\sum_{n=0}^{\infty}g_n\equiv\sum_{n=0}^{\infty}R_i(z_i)g_n=0$
almost surely as well.  Furthermore, if $z_i\in(1,2)$, then
\[
R_1(z_1)\dots R_m(z_m) {\cal G}
\]
is an analytic function in $z_1,\dots,z_m$.  Let $\Lambda_{n,m}{\cal G}$ be the terms
in the expansion of this analytic function which contain $z_1\dots z_n$ but no
other $z_i$.  Note that each $\Lambda_{n,m}{\cal G}$ equals 0 almost surely, and that
\[
\lim_{m\to\infty}\Lambda_{n,m}{\cal G} = y_1\dots y_n h_n(x_1,\dots,x_n) = 0
\]
almost surely.  In fact, the convergence to 0 follows from (\ref{9.10}),
since the operator $\Lambda_{n,m}$ removes more and more terms as 
$m$ increases.
This proves Lemma \ref{unique-expansion}.
\qed

\subsection{Proof of Uniqueness for $p\geq1$}
By Lemma~\ref{uniquely-determined} the following lemma gives the required uniqueness for the
 case $p\geq 1$.

\begin{lemma}
\label{unique-expansion-p-greater}
Let $p\geq 1$.
In (\ref{9.1}), the functions $g_n$ are uniquely determined.
\end{lemma}

\textbf{Proof.}
Fix $(t,x)$
and suppose that $u_1(t,x),u_2(t,x)$ are two solutions.  Subtracting,
we obtain
\begin{eqnarray*}
Y(t,x) &=& u_1(t,x)-u_2(t,x)  \\
&=& \sum_{n=0}^{\infty}g^{(1)}_n(t,x) - \sum_{n=0}^{\infty}g^{(2)}_n(t,x) \\
&=:& \sum_{n=0}^{\infty}\tilde{g}_n(t,x)
\end{eqnarray*}
where $g^{(i)}_n(t,x)$ corresponds to the solution $u_i(t,x)$.  Also, each
$\tilde{g}_n$ can be expressed as a stochastic integral
\[
\tilde{g}_n(t,x) = \int_{\widetilde{D}_n} f_n(t,x;x^{(1)},\dots,x^{(n)})
    L(dx^{(1)})\dots L(dx^{(n)}).
\]
Adding these terms, and using the notation from the beginning of
Section \ref{s9}, we find that
\[
Y(t,x) = f_0(t,x) + \int_{R(M(D),D)} \tilde{f}_1(t,x; y) L(dy)
\]
where
\begin{eqnarray*}
\lefteqn{ \tilde{f}_1(t,x; y) =  f_1(t,x; y)   } \\
&& + \sum_{n=2}^{\infty}
\int_{\tilde{R}_{n-1}(y)}f_n(t,x; x^{(1)},\dots,x^{(n-1)},y)L(dx^{(1)})\dots L(dx^{(n-1)})
\end{eqnarray*}
and $\tilde{R}_{n-1}(y)$ was defined in (\ref{8.5}).

Since $Y(t,x)=0$ and 
\[
z \to \int_{R(z,D)} \tilde{f}_1(t,x; y) L(dy)
\]
is a local $\mathcal{H}_z$ martingale, it follows that $f_0(t,x)=0$ and that the 
integrand $\tilde{f}_1(t,x; y)=0$ for almost every value of $y$.  But 
$\tilde{f}_1(t,x; y)=0$ is itself a sum of a deterministic function and a 
stochastic integral.  Therefore we may use the same argument to show that 
$f_1(t,x; y)=0$ for almost every value of $y$.  Continuing, we can use 
induction to show that for each value of $n$, $f_n(t,x; y_1,\dots,y_n)=0$ for 
almost every value of $y_1,\dots,y_n$.  
\qed

\appendix

\section{A Derivation of Hu's Estimate}
\setcounter{equation}{0}

In the proof of Theorem 4.1, equation (4.10) of \cite{hu02}, Hu gives the 
following estimate (in his equation, $\alpha=1-d/4$, and we use $\beta=d/4$).
\[
\mathcal{H}_n(t) 
:= \int_{\mathbf{T}_n(t)}\prod_{k=1}^{n}(s_{k+1}-s_k)^{-\beta}ds.
\]
Then
\begin{equation}
\label{10.1}
\mathcal{H}_n(t) = \frac{t^{n(1-\beta)}\Gamma(1-\beta)^n}{\Gamma(1+n(1-\beta))}.
\end{equation}
(Hu uses $\beta=\gamma/4$.  Also, Hu has 2 instead of 1 in the above formula.  
His answer is equivalent if we ignore multiplicative factors of n.)

Making the change of variables $s = \alpha r$, $ds = \alpha dr$,
we find
\begin{eqnarray*}
\mathcal{H}_n(\alpha t) &=& \int_{\mathbf{T}_n(\alpha t)}
  \prod_{k=1}^{n}(s_{k+1}-s_k)^{-\beta}ds\\
&=& \alpha^{-n\beta}\alpha^{n}\mathcal{H}_n(t)  \\
&=& \alpha^{n(1-\beta)}\mathcal{H}_n(t).
\end{eqnarray*}
Then we get the recurrence
\begin{eqnarray*}
\mathcal{H}_{n+1}(t) 
&=& \int_{0}^{t} (t-s)^{-\beta} \mathcal{H}_n(s)ds  \\
& = & \int_{0}^{t} (t-s)^{-\beta} s^{n(1 - \beta)}\mathcal{H}_n(1)ds  \\
& = & \mathcal{H}_n(1)\int_{0}^{t} (t-s)^{-\beta} s^{n(1 - \beta)}ds  \\
&=& t^{-n(1-\beta)}\mathcal{H}_n(t) \int_{0}^{t}(t-s)^{-\beta}s^{n(1-\beta)}ds.
\end{eqnarray*}
Let
\begin{eqnarray*}
a &=& \beta \\
b &=& n(1-\beta).
\end{eqnarray*}
Recall that the Euler beta function is defined as follows, for $a,b>0$.  
\[
B(a,b) := \int_{0}^{1}(1-x)^{a-1}x^{b-1}dx.
\]
It is a fundamental identity that
\[
B(a,b) = \frac{\Gamma(a)\Gamma(b)}{\Gamma(a+b)}.
\]
Changing variables, we get
\[
\int_{0}^{t}(t-s)^{-a}s^bds = \frac{t^{1-a+b}\Gamma(1-a)\Gamma(1+b)}
  {\Gamma(2-a+b)}.
\]
Therefore,
\begin{eqnarray*}
\mathcal{H}_{n+1}(t) &=& 
t^{-n(1-\beta)}\mathcal{H}_n(t) t^{1-\beta+n(1-\beta)}
  \frac{\Gamma(1-\beta)\Gamma(1+n(1-\beta))}{\Gamma(2-\beta+n(1-\beta))}   \\
&=& t^{(1-\beta)}\mathcal{H}_n(t)
  \frac{\Gamma(1-\beta)\Gamma(1+n(1-\beta))}{\Gamma(1+(n+1)(1-\beta))}   \\
\end{eqnarray*}
The product telescopes.  Since $\mathcal{H}_1(t) = t^{1-\beta}$, we get
(\ref{10.1}).

\end{document}